\documentclass{article}

\usepackage{t1enc}
\usepackage[latin2]{inputenc}

\usepackage{amssymb,amsmath}
\usepackage[slovene, croatian, slovak, english]{babel}

\newenvironment{proof}{\begin{trivlist}\item[]{\it
Proof.}}{\hfill$\square$\end{trivlist}}

\newtheorem{theorem}{Theorem}[section]
\newtheorem{corollary}[theorem]{Corollary}
\newtheorem{conj}[theorem]{Conjecture}

\newtheorem{lemma}[theorem]{Lemma}
\newtheorem{proposition}[theorem]{Proposition}

\newcommand{\C}{\mathbb C}
\newcommand{\la}{\lambda}

\newcommand{\haf}{\mathrm {haf}}

\newcommand{\per}{\mathrm {per}}

\begin{document}
\overfullrule=5pt

\date{}

\title{Remarks on the $\alpha$--permanent}

\author{ P\'eter E.\ Frenkel\footnote{Support by FNS and by OTKA grants  K 61116 and NK 72523 is gratefully acknowledged.}\\Universit\'e de Gen\`eve\\ 
2-4 rue du Li\`evre, 1211 Gen\`eve 4 \\{\tt frenkelp@renyi.hu}} 

\maketitle
\begin{abstract} 
We recall Vere-Jones's definition of the $\alpha$--permanent and describe the connection between the (1/2)--permanent and the hafnian. We establish expansion formulae for the $\alpha$--permanent in terms of partitions of the index set, and we use these to prove Lieb-type inequalities for the $\pm\alpha$--permanent of a positive semi-definite Hermitian $n\times n$ matrix and the $\alpha/2$--permanent of a positive semi-definite real symmetric $n\times n$ matrix if $\alpha$ is a nonnegative integer or $\alpha\ge n-1$. We are unable to settle Shirai's nonnegativity conjecture for  $\alpha$--permanents when $\alpha\ge 1$, but we verify it up to the $5\times 5$ case, in addition to recovering and refining some of Shirai's partial results by purely combinatorial proofs.


\end{abstract} 


\section{  Introduction}
Following Vere-Jones [V1, V2], we define the $\alpha$--permanent of the $
n\times
n$ matrix $A=(a_{i,j})\in M_n(\C)$ to be $$\per_\alpha A=\sum_{\pi\in\mathfrak S_n}\alpha ^{\nu(\pi)}\prod_{i=1}^na_{i,\pi(i)},$$
where $\mathfrak S_n$ is the symmetric group on $n$ elements and $\nu(\pi)$ is the number of disjoint cycles of the permutation $\pi$. In particular, $\alpha =1$ yields the ordinary permanent and $\alpha =-1$ yields $(-1)^n$ times the determinant. For real symmetric matrices, the case $\alpha=1/2$ recovers another known concept.  Recall that the hafnian  of a $2n\times
2n$ symmetric matrix $C=(c_{i,j})$ is  defined by $$\haf
\;C=\frac 1{n!2^n}\sum_{\pi\in\mathfrak S_{2n}} c_{\pi(1),\pi(2)}\cdots 
c_{\pi(2n-1),\pi(2n)}=\sum c_\Gamma,$$
where $\Gamma$ runs over the 1-regular graphs (perfect matchings) on $[2n]=\{1,\dots, 2n\}$ and \begin{equation*}c_\Gamma=\prod_{e\in E(\Gamma)}c_e;\end{equation*} note that we write $c_e=c_{ij}$ if $i$ and $j$ are the endpoints of the edge $e$.

\begin{proposition}\label{haf}  Let $A$ be a real symmetric $n\times n$ matrix. Then
\begin{equation}\label{Wick}
\per_{1/2} A=\frac{1}{2^n} \haf\left(\begin{matrix}
      A&A\\A&A\end{matrix}\right).
\end{equation}
\end{proposition}

This is essentially known. Since both sides are polynomials in the entries of $A$, we may assume that $A$ is positive semi-definite. Then we may consider centered, jointly Gaussian random variables $X_1$, \dots, $X_n$ with covariance matrix
$A$. The left hand side of~\eqref{Wick} 
is known to be equal to $2^{-n}E(X_1^2\cdot\cdots\cdot X_n^2)$, cf.\ Lu and Richards [LR] and Shirai [Sh]. The right hand side is the same quantity by the 
well-known [B, F, G, S, Z]
 Wick formula.

Nevertheless, a direct combinatorial proof may be of some interest.

\begin{proof}
Both sides are linear combinations of the $a_\Gamma$ where $\Gamma$ runs over the 2-regular graphs on the vertex set $[n]$.  We shall compute the coefficients of each side.

Suppose that $\Gamma$ has $\nu=\sum\nu_i$ connected components, of which $\nu_i$ are  $i$--cycles.  Then $n=\sum i\nu_i$.

The coefficient of $a_\Gamma$ on the LHS is $2^{-\nu}$ times the number of (1,1)--regular directed graphs whose underlying undirected graph is $\Gamma$. This number is \begin{equation*}
2^{-\nu}2^{\nu-\nu_1-\nu_2}=2^{-\nu_1-\nu_2},
\end{equation*}
 since each cycle of length $\ge 3$ has two non-isomorphic orientations, whereas each cycle of length 1 or 2  has only one.

The coefficient of $a_\Gamma$ on the RHS is $2^{-n}$ times the number of 1--regular  graphs on $[2n]\simeq [n]\times [2]$ that project to $\Gamma$. A cycle of length 1 (i.e., a loop) in $\Gamma$ has only one lifting to $[n]\times [2]$. A cycle of length 2 (i.e., two parallel edges) has two liftings.  A cycle of length $i\ge 3$ has $2^i$ liftings. Thus, the coefficient is
\begin{equation*} 2^{-n}2^{\nu_2}\prod_{i\ge 3}2^{i\nu_i}=2^{-n+\nu_2+n-\nu_1-2\nu_2}=2^{-\nu_1-\nu_2},
\end{equation*} and the Proposition follows.
\end{proof}

\section{Expansion formulae}\label{exp}
In this section, we shall expand $\per_\alpha A$ in terms of certain $\beta$--permanents of diagonal submatrices of $A$. These expansions shall be described in terms of partitions of the set $[n]$. When $A$ is the identity matrix, our arguments reduce more or less to some classical ideas of Gian-Carlo Rota related to enumerating set partitions [R].

Put $[n]=\{1,\dots, n\}$. 
For an $n\times n$ matrix $A=(a_{i,j})$ and a subset $I$ of $[n]$, we write $A[I]:=(a_{i,j})_{i,j\in
I}.$ The symmetric group on $I$ is written $\mathfrak S(I)$.
\begin{lemma} We have
\begin{equation}\label{sum}
\per_{\beta_1+\dots+\beta_m} A=\sum\prod_{j=1}^m\per_{\beta_j}A[I_j],
\end{equation}
the summation being over all ordered partitions $(I_1,\dots, I_m)$ of $[n]$ into $m$ disjoint (possibly empty) subsets.
\end{lemma}

\begin{proof} For any permutation $\pi$, let us write $\Pi(\pi)=\{C_1,\dots, C_{\nu(\pi)}\}$ for the unordered partition given by the cycles of $\pi$ (these are non-empty subsets of $[n]$).  Then
\begin{align*} \per_{\sum\beta_j} A=\sum_\pi\left(\sum\beta_j\right)^{\nu(\pi)}\prod_ia_{i,\pi(i)}=\\=\sum_\pi\prod_{C\in\Pi(\pi)}\left(\sum\beta_j\prod_{i\in C}a_{i,\pi(i)}\right)=\\=\sum_\pi\sum_{f:\Pi(\pi)\to [m]}\prod_{C\in\Pi(\pi)}\left(\beta_{f(C)}\prod_{i\in C}a_{i,\pi(i)}\right)=\\=
\sum_{f:[n]\to [m]}\sum_{f\circ \pi =f}\prod_{C\in\Pi(\pi)}\left(\beta_{f(C)}\prod_{i\in C}a_{i,\pi(i)}\right)=\\=\sum_{(I_1,\dots, I_m)}\sum_{\pi_1\in\mathfrak S(I_1)}
\cdots\sum_{\pi_m\in\mathfrak S(I_m)} \prod_{j=1}^m \left(\beta_j^{\nu(\pi_j)}\prod_{i\in I_j}a_{i,\pi(i)}\right)=\\=\sum\prod_{j=1}^m\per_{\beta_j}A[I_j].
\end{align*}
\end{proof}

We shall now apply the lemma to the case where $\beta_1=\dots=\beta_m$.  It will be convenient to get rid of the empty subsets appearing in the partitions.

Let us define
\begin{equation}
\per_\beta (A,k)=\sum_{(I_1,\dots, I_k)}\prod_{j=1}^k\per_\beta A[I_j],
\end{equation} the summation being over all oredered partitions
$(I_1,\dots, I_k)$ of $[n]$ into $k$ disjoint, \emph{nonempty} subsets. We abbreviate $\per_1$ to per and $\per_{-1}$ to $(-1)^n\det$.

As usual, we define $\binom\alpha k=\alpha(\alpha-1)\cdots (\alpha-k+1)/k!$.

\begin{theorem} For any numbers $\alpha$ and $\beta$, and any $n\times n$ matrix $A$, we have
\begin{equation}\label{product}
\per_{\alpha\beta} A=\sum_{k=1}^n\binom\alpha k\per_\beta (A,k).
\end{equation}

In particular,
\begin{equation}\label{pos}
\per_{\alpha} A=\sum_{k=1}^n\binom\alpha k\per (A,k)
\end{equation}

and 

\begin{equation}\label{neg}
\per_{-\alpha} A=(-1)^n\sum_{k=1}^n\binom\alpha k\det (A,k).
\end{equation}
Also, if $A$ is real and symmetric,

\begin{equation}\label{half}
\begin{aligned}
\per_{\alpha/2} A =\\=
\frac 1{2^n}\sum_{k=1}^n\binom\alpha k\sum
\left\{
\prod_{j=1}^k\haf\left(\begin{matrix}A[I_j]&A[I_j]\\A[I_j]&A[I_j]\end{matrix}\right)
\mid
\coprod_{j=1}^k I_j=[n], \forall I_j\ne\emptyset\right\}
.
\end{aligned}
\end{equation}
\end{theorem}

\begin{proof}
Both sides  are polynomials in $\alpha$, so we may assume that $\alpha=m$ is a nonnegative integer. By the Lemma, we have \begin{align*}
\per_{\alpha\beta} A=\per_{\beta+\dots+\beta} A =\sum\prod_{j=1}^\alpha\per_\beta A[I_j],\end{align*}
where we are summing over ordered partitions with  empty subsets allowed.  However, the $\beta$--permanent of the empty matrix is 1. Thus we may restrict ourselves to ordered partitions with nonempty subsets, and such a partition, if it has $k$ parts, will be obtained $\binom\alpha k$ times. Hence the result.
\end{proof}

\section{Inequalities for positive semi-definite matrices}\label{ineq}
Throughout this section, $A$ will be a positive semi-definite Hermitian $n\times n$ matrix. 

Shirai and Takahashi [Sh, ShT]
have conjectured that $\per_\alpha A\ge 0$ if $\alpha\ge 1$, and that $\per_{\alpha/2}A\ge 0$ if $\alpha\ge 1$
 and $A$ is real. 
Shirai  [Sh] proves $\per_\alpha A\ge 0$ if $\alpha$ is a nonnegative integer or $\alpha\ge \mathrm{rank}( A)-1$,  proves $\per_{\alpha/2} A\ge 0$ for real $A$ if $\alpha$ is a nonnegative integer or $\alpha\ge n-1$, and    proves
also that $(-1)^n\per_{-\alpha}A\ge 0$ if $\alpha$ is a nonnegative integer.

The question of nonnegativity is motivated by problems from probability theory. See [EK, Sh, ShT, V1, V2] and references therein.
Note that [EK, Sh, ShT] formulate everything in the terms of the $\alpha$--determinant
\[\mathrm{det}_\alpha A=\alpha^n\per_{1/\alpha}A\]
rather than the $\alpha$--permanent used by Vere-Jones and
the present paper.

We shall now strengthen some of Shirai's nonnegativity results to obtain Lieb  type inequalities when  $\alpha$ is an integer or $\alpha\ge n-1$. 
 Also, we 
verify that $\per_\alpha A\ge 0$ if $\alpha\ge 1$ and $n\le 5$. Sadly, the conjectures of Shirai and Takahashi remain open in general. Nevertheless, we propose a stronger conjecture.

Suppose that the p.s.d.H. matrix $A$ is partitioned as
\begin{equation}\label{bl}A=\left(\begin{matrix}
      A'&B\\B^*&A''\end{matrix}\right).\end{equation}
Put
\begin{equation}D=\left(\begin{matrix}
      A'&0\\0&A''\end{matrix}\right).\end{equation}
Recall Lieb's inequality [L, D, Mi]
\begin{equation}\label{Lieb} 
\per A\ge\per D
=\per A'\cdot \per
A''\end{equation}
and the classical Fischer inequality
\begin{equation}\label{Fischer} 
\det A\le\det D
=\det A'\cdot \det
A''.\end{equation}
We immediately deduce
\begin{equation}\label{L} 
\per (A,k)\ge\per (D,k)\end{equation}
and 
\begin{equation}\label{F} 
\det (A,k)\le\det (D,k).\end{equation}

When $A$ is real, recall from [F] the the inequality
\begin{equation}\label{Frenkel} 
\haf \left(\begin{matrix}
      A&A\\A&A\end{matrix}\right)\ge\per A.\end{equation}
We deduce
\begin{equation}\label{Fr}
\sum
\left\{
\prod_{j=1}^k\haf\left(\begin{matrix}A[I_j]&A[I_j]\\A[I_j]&A[I_j]\end{matrix}\right)
\mid
\coprod_{j=1}^k I_j=[n], \forall I_j\ne\emptyset\right\}\ge\per (A,k).
\end{equation}

\begin{theorem}\label{Lieb-type}
Suppose that $\alpha$ is a nonnegative integer or $\alpha\ge n-1$. Then, for $A\in M_n(\C)$ p.s.d.H.
partitioned as in \eqref{bl},
we have
\begin{equation}
\per_\alpha A\ge\per_\alpha D=\per_\alpha A'\cdot\per_\alpha A''\end{equation}
and 
\begin{equation}\label{csillagos}
0
\le
(-1)^n\per_{-\alpha} A\le(-1)^n\per_{-\alpha} D=(-1)^n\per_{-\alpha} A'\cdot\per_{-\alpha} A''.\end{equation}
If, in addition, $A$ is real, then
\begin{equation}\label{real}
\per_{\alpha/2}A\ge\frac1{2^n}\per_\alpha A.
\end{equation}
\end{theorem}

\begin{proof}
The assumption on $\alpha$ ensures that $\binom\alpha k\ge 0$ for $k=1, \dots, n$. Thus, by \eqref{pos} and \eqref{L},
\begin{align*}\per_\alpha A&=\sum_{k=1}^n\binom\alpha k\per (A,k)\ge\\&\ge\sum_{k=1}^n\binom\alpha k\per (D,k)=\per_\alpha D=\per_\alpha A'\cdot\per_\alpha A''.\end{align*}

Similarly, by \eqref{neg} and \eqref{F},
\begin{align*}0\le (-1)^n\per_{-\alpha} A=\sum_{k=1}^n\binom\alpha k\det (A,k)\le\\\le\sum_{k=1}^n\binom\alpha k\det (D,k)=(-1)^n\per_{-\alpha} D=(-1)^n\per_{-\alpha} A'\cdot\per_{-\alpha} A''.\end{align*}

Finally, if $A$ is real, then by \eqref{half} and \eqref{Fr},
\begin{align*}
\per_{\alpha/2} A =\\=
\frac 1{2^n}\sum_{k=1}^n\binom\alpha k\sum
\left\{
\prod_{j=1}^k\haf\left(\begin{matrix}A[I_j]&A[I_j]\\A[I_j]&A[I_j]\end{matrix}\right)
\mid
\coprod_{j=1}^k I_j=[n], \forall I_j\ne\emptyset\right\}\ge\\\ge\frac 1{2^n}\sum_{k=1}^n\binom\alpha k\per (A,k)=\frac 1{2^n}\per_\alpha A.
\end{align*}
\end{proof}

\begin{corollary}
Suppose that $\alpha$ is  a  nonnegative integer or $\alpha\ge n-1$. Then, for $A\in M_n(\C)$ p.s.d.H., we have
\begin{equation}\label{Marcus}
\per_\alpha A\ge\alpha ^n\prod_{i=1}^na_{ii}\ge(-1)^n\per_{-\alpha} A
.
\end{equation}
If, in addition, $A$ is real, then
\begin{equation}\label{half-Marcus}
\per_{\alpha/2} A\ge (\alpha/2)^n\prod_{i=1}^na_{ii}.
\end{equation}
\end{corollary}

\begin{proof}
Obvious induction to prove \eqref{Marcus}, then \eqref{real} and \eqref{Marcus} to prove \eqref{half-Marcus}.
\end{proof}

\begin{conj}
The condition on $\alpha$ can be relaxed to $\alpha\ge 1$ for all inequalities stated in Theorem~\ref{Lieb-type} and its Corollary, except for the 
leftmost inequality in \eqref{csillagos}.
\end{conj}

To support the conjecture, we prove the inequalities~\eqref{Marcus} of the above corollary for small matrices under the relaxed condition for $\alpha$.

\begin{theorem}
Suppose that $\alpha\ge 1$ and $n\le 5$. Then, for $A\in M_n(\C)$ p.s.d.H., the inequalities~\eqref{Marcus} hold.
\end{theorem}

\begin{proof} We may assume that $a_{ii}=1$ for all $i$.

If the Theorem is true for a number $\alpha$, then by the Lemma, it is also true for $\alpha+1$. We may therefore assume that $1\le\alpha\le 2$.

We may assume that $n=5$ since the statement for $A$ is equivalent to that for $A\oplus\mathbf 1_{5-n}$.

From now on $\bf 1$ is the $5\times 5$ identity matrix.

The statement to be proven is $\per_{\pm\alpha }A\ge\per_{\pm\alpha }\mathbf 1$.

In view of formula~\eqref{product}, it suffices to prove that
\begin{equation}\label{123}
\sum_{k=1}^3\binom\alpha k\per_{\pm 1} (A,k)\ge\sum_{k=1}^3\binom\alpha k\per_{\pm 1} (\mathbf 1,k)
\end{equation}
and
\begin{equation}\label{45}
\sum_{k=4}^5\binom\alpha k\per_{\pm 1} (A,k)\ge\sum_{k=4}^5\binom\alpha k\per_{\pm 1} ({\mathbf {1}},k).
\end{equation}

For a partition $\lambda=(\la_1\ge\la_2\ge\dots\ge\la_k) $ of 5 into $k$ positive integer parts, we define $p(\la)$ to be the average value of
\begin{equation*}
\prod_{j=1}^k\per_{\pm 1} A[I_j],
\end{equation*}
where we are averaging over the partitions $(I_1,\dots, I_k)$ of $[5]$ such that $|I_j|=\la_j$ for $j=1,\dots, k$.
From the Lieb and Fischer inequalities, 
we get $p(\la)\ge p(\mu)$ if $\la$ arises from $\mu$ by replacing two parts of $\mu$ by their sum.

Then \eqref{45} reduces to
\begin{equation*}
 10p(2,1,1,1)
+(\alpha -4)p(1,1,1,1,1)\ge\pm(10+\alpha-4),\end{equation*}
which is true because \begin{equation*}
p(2,1,1,1)\ge
p(1,1,1,1,1)=\pm 1. \end{equation*}
 
Also, \eqref{123} reduces to
\begin{equation}\label{123red}
 \begin{aligned}
p(5)+(\alpha-1)(5p(4,1)+10p(3,2))+\\
+(\alpha-1)(\alpha-2)(10p(3,1,1)+15p(2,2,1))
\ge\\\ge
\pm(1+(\alpha-1)(5+10)+(\alpha-1)(\alpha-2)(10+15
)).
\end{aligned}
\end{equation}
By the Lieb and Fischer inequalities, we have $$p(2,2,1)\le p(3,2)$$
and $$p(3,1,1)\le\min(p(4,1),p(3,2))\le\frac56 p(4,1)+\frac 16p(3,2).$$
Thus, the LHS of \eqref{123red} is at least
\begin{equation}\label{LHS}
p(5)+5(\alpha-1)\left(1+\frac53(\alpha-2)\right)(p(4,1)+2p(3,2)).
\end{equation}
Now
$$p(5)\ge\max(p(4,1),p(3,2))\ge\frac13(p(4,1)+2p(3,2)),$$
so \eqref{LHS} is at least
\begin{align*}
\left(\frac13 +5(\alpha-1)\left(1+\frac53(\alpha-2)\right)\right)(p(4,1)+2p(3,2)).
\end{align*}
Here the first factor is non-negative and the last factor is at least $\pm 3$, whence the result.
\end{proof}


\section*{References}
\noindent
[B] A.\ Barvinok, Integration and optimization of multivariate polynomials by
restriction onto a random subspace. Found.\ Comput.\ Math.\ 7 (2007) no 2., 229--244.

\bigskip\noindent
[D] D.\ \v Z.\ \DJ okovi\'c, Simple proof of a
theorem on permanents, Glasgow
Math.\ J.\ 10 (1969), 52--54.

\bigskip\noindent
[EK] N.\ Eisenbaum,  H.\ Kaspi,
On permanental processes. Stochastic Process.\ Appl.\ 119 (2009), no. 5, 1401--1415. 

\bigskip\noindent
[F] P.\ E.\ Frenkel, Pfaffians, Hafnians and products of real linear functionals.  Math. Res. Lett.  15  (2008),  no. 2, 351--358. 

\bigskip\noindent
[G] L.\ Gurvits, Classical complexity and quantum entanglement.  J.\ Comput.\
 System Sci.\  69  (2004),  no. 3, 448--484.

\bigskip\noindent
[L] E.\ H.\ Lieb, Proofs of some conjectures on permanents, J.\ Math.\ Mech.\
16 (1966), 127--134.
\bigskip\noindent
[LR] I.\ L.\ Lu and D.\ P.\ Richards, MacMahon's master theorem, representation theory, and moments of Wishart distributions, Adv.\ Appl.\ Math.\ 27(2001), 531--547.

\bigskip\noindent
[Mi] H.\ Minc, Permanents, Encyclopedia of Mathematics and its Applications,
Add\-is\-on-Wesley, 1978

\bigskip\noindent
[R] G.-C.\ Rota, The number of partitions of a set, Amer.\ Math.\ Monthly 71, no.\ 5 (May 1964), 498--504.

\bigskip\noindent
[Sh] T.\ Shirai, Remarks on the positivity of alpha-determinants, Kyushu J. Math. 61 (1) (2007),  169--189.

\bigskip\noindent
[ShT] T.\ Shirai and Y.\ Takahashi, Random point fields associated with certain Fredholm determinants I: Fermion, Poisson and boson point processes, J. Funct. Anal. 205 (2003),  414--463. 

\bigskip\noindent
[S] B.\ Simon,
The P$(\phi)_2$ Euclidean (Quantum) Field Theory, Princeton Series in Physics,
Princeton University Press, 1974

\bigskip\noindent
 [V1] D.\ Vere-Jones, A generalization of permanents and determinants, Linear Alg.\ Appl.\ 111 (1998), 119--124.

\bigskip\noindent
 [V2] D.\ Vere-Jones, Alpha-permanents and their applications to multivariate Gamma, negative binomial and ordinary binomial distribution. New Zealand J.\ of Math.\ 26 (1997), 125--149.

\bigskip\noindent
[Z] A.\ Zvonkin, Matrix integrals and map enumeration: an accesible
introduction, Combinatorics and physics (Marseille, 1995), Math.\ Comput.\
Modelling 26 (1997), 281--304.
\end{document}